\documentclass[reqno]{amsart}
\usepackage[LGR,T1]{fontenc}

\usepackage{comment,todonotes,ifthen}
\presetkeys{todonotes}{inline}{}

\usepackage{epigraph}
\DeclareMathOperator{\Sym}{Sym}

\newcommand{\SL}{\textup{SL}}
\newcommand{\GL}{\textup{GL}}

\newcommand{\C}{\mathbb{C}}
\newcommand{\F}{\mathbb{F}}

\DeclareMathOperator{\Par}{Par}

\renewcommand{\emptyset}{\varnothing}
\usepackage{bm}

\usepackage{Style}
\usepackage{ytableau}
 
\newcommand{\qbinom}[2]{{\genfrac{[}{]}{0pt}{}{#1}{#2}}}
\usepackage{mathtools}

\makeatletter
\renewenvironment{smallmatrix}[1][.2em]{\null\,\vcenter\bgroup
  \Let@\restore@math@cr\default@tag
  \baselineskip6\ex@ \lineskip1.5\ex@ \lineskiplimit\lineskip
  \ialign\bgroup\hfil$\m@th\scriptstyle##$\hfil&&\kern#1\hfil
  $\m@th\scriptstyle##$\hfil\crcr
}{%
  \crcr\egroup\egroup\,%
}
\makeatother

\newcommand{\showcommentsbox}{yes}
\newsavebox{\commentbox}
\newenvironment{com}%
{\ifthenelse{\equal{\showcommentsbox}{yes}}%
    {\footnotemark
        \begin{lrbox}{\commentbox}
            \begin{minipage}[t]{1.25in}\raggedright\sffamily\tiny
                \footnotemark[\arabic{footnote}]}
                {\begin{lrbox}{\commentbox}}}%
                {\ifthenelse{\equal{\showcommentsbox}{yes}}%
                {\end{minipage}\end{lrbox}\marginpar{\usebox{\commentbox}}}
        {\end{lrbox}}}

\renewcommand{\le}{\leqslant}
\renewcommand{\leq}{\leqslant}
\renewcommand{\ge}{\geqslant}
\renewcommand{\geq}{\geqslant}

\title[On the generalised Foulkes conjecture]{On the generalised Foulkes conjecture\\ for $\mathrm{SL}_2(\mathbb{C})$ under divisibility conditions}

\author[M.~Gangl \and \'A.~Guti\'errez \and M.~Szwej]{Moritz Gangl \and \'Alvaro Guti\'errez \and Micha{\l} Szwej}
\address{University of Vienna}
\email[M.~Gangl]{moritz.gangl@univie.ac.at}
\address{University of Bristol}
\email[\'A.~Guti\'errez]{a.gutierrezcaceres@bristol.ac.uk}
\email[M.~Szwej]{michal.szwej@bristol.ac.uk}

\date{\today}

\begin{document}

\begin{abstract}
    The generalised Foulkes conjecture for $\mathrm{SL}_2(\mathbb{C})$ was recently proved under mild divisibility conditions by Raicu, Sam, Weyman, and Yang, who showed surjectivity of the Foulkes--Howe map through a geometric and homological approach. As an immediate corollary, we note partial proofs of related conjectures by Bergeron, Zanello, and Troyka. The main result of this paper is that the dual of the (geometric) Foulkes--Howe map is the (combinatorial) $k$-fold plethystic substitution, which admits an explicit and straightforward definition. We derive several formulas and combinatorial interpretations for its structure constants. We briefly remark on an unexpected corollary that settles a conjecture in condensed matter physics. Finally, we use the combinatorial properties of the $k$-fold map to propose candidate maps towards other variants of Foulkes' conjecture.\medskip
    
\end{abstract}
\maketitle
\vspace{-.5cm}
\section{Introduction}
Consider $\SL_N(\C)$ acting in its natural representation on $\C^N$.
The irreducible representations of $\SL_N(\C)$ can be constructed as $S^\lambda \C^N$, where $S^\lambda$ is a \emph{Schur functor} and $\lambda$ ranges over partitions of length strictly less than $N$. A major open problem in representation theory is the \emph{plethysm problem} \cite{StanleyList, Mystery}: how does $S^\nu S^\mu \C^N$ decompose into irreducible representations? More precisely, the goal is to understand the coefficients $a^\lambda_{\nu[\mu]}$ occurring in the decomposition
\[
S^\nu S^\mu \C^N = \bigoplus\nolimits_{\lambda} \big(S^\lambda\C^N\big)^{\oplus a^\lambda_{\nu[\mu]}}.
\]

When $\lambda = (n)$ is a one-row partition, the Schur functor $S^{(n)}$ is the more familiar symmetric power functor $\Sym^n$.
Although much simpler than the general problem, even the following conjecture remains open \cite{Foulkes, McKay, EvseevPagetWildon, PW2, CIM}. Throughout this paper we write $U\leq V$ when $U$ is isomorphic to a subrepresentation of~$V$.
\begin{conjecture}[Foulkes]\label{conj:0}
    If $m \le n$ then $\Sym^m \Sym^n \C^N \le \Sym^n \Sym^m \C^N$.
\end{conjecture}
Several generalisations of Foulkes' conjecture are proposed by Doran \cite{Doran}, Ab\-des\-se\-lam and Chipalkatti \cite{AC,AC14}, and 
Bergeron \cite{Bergeron}.
A four parameter generalisation due to Vessenes \cite{Vessenes} is of special interest to this paper.
\begin{conjecture}[Generalised Foulkes conjecture] \label{conj:A}
    Let $a, b, c, d \in \N$, such that $ab=cd$, and assume $a=\min(a,b,c,d)$. Then
    \(
    \Sym^b\Sym^a\C^N \le \Sym^d\Sym^c\C^N.
    \)
\end{conjecture}

For $\SL_2(\C)$, Foulkes' conjecture holds (with equality) by Hermite reciprocity:
\[
\Sym^m \Sym^n \C^2 \cong \Sym^n \Sym^m \C^2,
\]
but the generalised Foulkes conjecture is highly non-trivial even in this case.
It is known to hold for $a = 2$ \cite{AC} and $a = 3$ \cite{Vessenes, Zanello}. 
In a recent preprint, Raicu, Sam, Weyman, and Yang proved the generalised Foulkes conjecture for $\SL_2(\C)$ under mild divisibility conditions.
\begin{thm}[{\cite[Theorem 1.4]{RSWY}}]
\label{thm:main}
The \emph{Foulkes--Howe map} 
$$\alpha_d: \Sym^d\Sym^{ka}\C^2\rightarrow\Sym^{kd}\Sym^a\C^2$$
is injective if $d\leq a$ and surjective if $d\geq a$.
In particular, Conjecture~\ref{conj:A} holds for $\SL_2(\C)$ when $a$ divides either $c$ or $d$.
\end{thm}
The Foulkes--Howe map was constructed by Brion \cite{Brion} in his study of Conjecture~\ref{conj:0} as a morphism of algebraic group representations, by Ab\-des\-se\-lam and Chipalkatti~\cite{AC} in the study of Conjecture~\ref{conj:A} for $\SL_2(\C)$ as an explicit manipulation on the coordinates of certain binary forms, and by Raicu, Sam, Weyman, and Yang \cite{RSWY} as a map between sheaf cohomology groups in the study of Hilbert covariants (see also \cite{AC14} where this approach was pioneered).

Another candidate towards Conjecture~\ref{conj:A}, the $k$-\emph{fold map}, going in the opposite direction to $\alpha_d$, was considered combinatorially in \cite{GMSW-fpsac}. It relies on the realisation of plethysm of symmetric powers 
\[
\Lambda^{\F}_{\le n}[\mathbf{x}_m] \cong \Sym_m \Sym^n \F^2
\]
from \cite{GMSW}, where $\F$ is a field of arbitrary characteristic and $\Lambda^{\F}_{\leq n}[\mathbf{x}_m]$ is the $\F$-vector space of symmetric polynomials of degree at most $n$ in each of the $m$ variables. In this identification, the $k$-fold map is defined explicitly as 
\begin{align}
    \kappa_{a,d}^k:\Lambda^{\F}_{\le a}[\mathbf{x}_{kd}]&\rightarrow\Lambda^{\F}_{\le ka}[\mathbf{y}_d]\\
    f(\mathbf{x}_{kd})& \mapsto f(y_1,\overset{k}{\ldots}\,, y_1,\ldots,y_d,\overset{k}{\ldots}\,,y_d).
\end{align}
Previously, this map had been considered in the study of the particle entanglement spectrum of Laughlin states \cite{GEA}. More literature is available for the map $\kappa_{\infty,\infty}^k$ on symmetric functions \cite{Baker, Brenti,FJK}. It was pointed out by Abdelmalek Abdesselam and by Mark Wildon in private communication that the Foulkes--Howe map and the $k$-fold map might be related.\medskip

In this paper, we give a characteristic-independent definition of the Foulkes--Howe map $\alpha_d$, which over $\mathbb{C}$ agrees with the map constructed in \cite{RSWY}. Moreover, we show the following result over an arbitrary field $\F$.
\begin{thm}\label{duality}
    Let $a,d,k\in\mathbb{N}$ with $a\leq d$. The Foulkes--Howe map
    $$\begin{tikzcd}
         \Sym^d\Sym_{ka}(\F^2)^{\star\star}\arrow[rr,"\alpha_d"]&&\Sym^{kd}\Sym_a(\F^2)^{\star\star}
     \end{tikzcd}$$
    is dual to the $k$-fold map
    \[
    \begin{tikzcd}[row sep=small]
    \Sym_{kd}\Sym^{a}(\F^2)^{\star} \arrow[rr, "\kappa_{a,d}^k"]                                        &  & \Sym_d\Sym^{ka}(\F^2)^{\star}                                          \\
    {\Lambda^{\F}_{\le a}[\mathbf{x}_{kd}]} \arrow[rr, "\kappa_{a,d}^k"'] \arrow[u, Rightarrow, no head] &  & {\Lambda^{\F}_{\le ka}[\mathbf{x}_{d}].} \arrow[u, Rightarrow, no head]
    \end{tikzcd}
    \]
\end{thm}
\begin{cor}[{\cite[Conjecture~7.3]{GMSW-fpsac}}]
  If $\F=\CC$ then the $k$-fold map
  $\kappa_{a,d}^k$ is injective for $a\le d$ and surjective for $a \ge d$.
\end{cor}
The proof of Theorem~\ref{duality} relies on a careful analysis of the dualities involved in the definition of both the Foulkes--Howe map and the $k$-fold map, and hence why we distinguish between a vector space and its double dual. The Foulkes--Howe map is ultimately a map between the spaces of twisted global sections of certain projective spaces of binary forms. We compute the coefficients of this map on a basis and recover the dual of a formula for the $k$-fold map. 
Whereas Schur and Weyl functors are dual and isomorphic functors over $\CC$, distinguishing between them (by working over an arbitrary field $\F$) is an increasingly successful tool to understand $\SL_2$-plethysms \cite{McDowellWildon,GMSW,IOT,GMSW-fpsac}. 

The $k$-fold map has many rich combinatorial properties, and its structure constants are of independent interest. In \S\ref{sec:matrix} we study the Schur and monomial expansions of the $k$-fold map of symmetric functions.
We exploit formulas found in~\cite{Baker, Brenti} as well as some new formulas to show that the Schur-basis matrix of the $k$-fold map $\kappa^k_{\infty,\infty}$ has many symmetries, whereas the monomial-basis matrix is lower triangular. We give Cauchy product formulas, as well as a combinatorial interpretation for the coefficients as stated in the following proposition; see \S\ref{sec:matrix} or \cite[\S3]{SXP} for the necessary definitions.
\begin{prop}\label{thm:Schur count}
    Let $\lambda,\mu$ be two partitions and $k\in\mathbb{N}$. The coefficient of $s_\mu$ in the Schur expansion of $\kappa^k_{\infty,\infty}(s_\lambda)$ is given by the number of pairs of lattice $k$-multi-tableaux with contents $\lambda$ and $\mu$ of an arbitrary $k$-multi-shape.
\end{prop}

Since $\SL_2(\CC)$ is a subgroup of $\GL_2(\CC)$, any $\GL_2(\CC)$ action restricts to an $\SL_2(\CC)$ action.
The study of $\GL_2(\CC)$ and $\SL_2(\CC)$ plethysms is intimately related to the study of $q$-analogues.
The $\GL_2(\CC)$-character of $\Sym^m\Sym^n\C^2$ is the principal specialisation
\[
s_{(m)}(1,q,\ldots,q^n) = \qbinom{n+m}{n}
\]
of the Schur function $s_{(m)}$ indexed by the partition $(m)$, where $\qbinom{n+m}{n}$ denotes the \emph{$q$-binomial coefficient} or the \emph{Gaussian coefficient}. 
For a precise definition, see \S\ref{sec:preliminaries}.
The study of $q$-analogues is a rich and growing field, and Gaussian coefficients are present in countless $q$-identities~\cite{Koepf} and
can be interpreted combinatorially as the generating functions for subsets of~$\N_0$ by the sum of their elements, partitions in an $n\times m$ rectangle, and inversions in permutations of $n$ objects
\cite[Ch.~1]{StanleyEC2}.

An immediate consequence of Theorem \ref{thm:main} is a partial proof of a conjecture of Bergeron \cite{BergeronSlides} whenever $a$ divides either $c$ or $d$.
\begin{cor}[$q$-Foulkes' conjecture for Gaussian coefficients] \label{cor: Bergeron}
Let $a, b,$ $c, d \in \N$, such that $ab=cd$, and assume $a=\min(a,b,c,d)$. 
If $a$ divides either $c$ or $d$, then $\qbinom{c+d}{c} - \qbinom{a+b}{a} \in \N_0[q]$.
\end{cor}
A polynomial $a_0 + a_1 q + \cdots + a_nq^n \in\N_0[q]$ is \emph{unimodal} if $a_1 \le \cdots \le a_k \ge \cdots\ge a_n$ for some $k$. Gaussian coefficients are  unimodal polynomials, which was established by O'Hara in her influential paper~\cite{OHara}.
The study of unimodal polynomials has recently grown in popularity due to the major successes of Huh and his collaborators~\cite{Huh}. Another immediate consequence of Theorem~\ref{thm:main} is a partial proof of a conjecture of Zanello~\cite{Zanello} under the usual divisibility conditions.
\begin{cor} \label{cor: Zanello}
Let $a, b, c, d \in \N$, such that $ab=cd$ and $a=\min(a,b,c,d)$. If $a$ divides either $c$ or $d$, then $\qbinom{c+d}{c} - \qbinom{a+b}{a}$ is a unimodal polynomial.
\end{cor}

In \S\ref{sec:Troyka} we derive a sufficient condition for $q$-positivity of $\qbinom{c+d}{c}-\qbinom{a+b}{b}$, for general $a,b,c,d$. As a corollary, we settle Troyka's conjecture (Conjecture~\ref{conj:Troyka}) under the usual divisibility conditions, vastly generalising Corollary~\ref{cor: Bergeron} by relaxing the equality $ab=cd$ to a necessary condition $ab\le cd$.
\begin{cor}\label{cor:troyka-div}
    Let $a,b,c,d\in\N$, such that $ab\le cd$ and $a=\min(a,b,c,d)$. If $a$ divides either $c$ or $d$ then $\qbinom{c+d}{c}-\qbinom{a+b}{a}\in\N_0[q]$.
  \end{cor}
  We further strengthen this result in Proposition~\ref{prop:BestTroyka}.

\bigskip

\textbf{Outline of the paper:} In Section~\ref{sec:preliminaries} we cover the basics of symmetric functions, $q$-analogues, and representations as vector spaces of polynomials necessary for the later parts of the paper. Section~\ref{sec:Troyka} discusses implications of Theorem~\ref{thm:main} for the conjectures related to Gaussian binomial coefficients. In Section~\ref{sec:main}, we introduce the \emph{$k$-fold map}, present some of its basic properties, and give a brief overview of its appearances in the literature. This is followed in Section~\ref{sec: proof of duality} by a recap of the Foulkes--Howe map from~\cite{RSWY} and the proof of Theorem~\ref{duality}. In Section~\ref{sec:matrix} we present our results on the structure constants of the $k$-fold map. We conclude the paper in Section~\ref{further} by proposing two extensions of the $k$-fold map: first, towards the full version of Conjecture~\ref{conj:A} (without the divisibility conditions), and second, towards Foulkes' conjecture (Conjecture~\ref{conj:0}; when $m$ divides $n$).

\section{Preliminaries}\label{sec:preliminaries}
\subsection{Symmetric functions}
We follow \cite[Ch.~7]{StanleyEC2} for background on symmetric functions, with the distinction that we work over $\C$ rather than $\ZZ$ or $\Q$.

An \emph{alphabet} is a set $\mathbf{x}_n = \{x_1, \ldots, x_n\}$ of variables.
Let
\(
\Lambda[\mathbf{x}_n] = \C[x_1, \ldots, x_n]^{S_n}
\)
be the \emph{$\C$-algebra of symmetric polynomials in $\mathbf{x}_n$}, which are polynomials in $x_1, \ldots, x_n$ invariant under any permutation of variables. 
Given an infinite alphabet $\mathbf{x}_\infty = \{x_1, x_2, \ldots\}$, the \emph{$\C$-algebra of symmetric functions in $\mathbf{x}_\infty$} is the graded colimit (or direct limit) as $n\to\infty$ of the $\C$-algebras of symmetric polynomials in $\mathbf{x}_n$.

A \emph{partition} $\lambda$ is a weakly decreasing finite sequence of non-negative integers $(\lambda_1, \lambda_2, \ldots, \lambda_k)$, called the \emph{parts} of $\lambda$. Its \emph{length} $\ell(\lambda)$ is the number of non-zero parts, and its \emph{size} is $|\lambda|=\lambda_1+\cdots+\lambda_k$. We assume $\lambda_i=0$ for any $i>\ell(\lambda)$.
We write $\Par$ for the set of partitions, $\Par_{\le h}$ for $\{\lambda\in\Par : \ell(\lambda)\le h\}$, and $L(h,w)$ for $\{\lambda\in\Par_{\le h} : \lambda_1\le w\}$.
The \emph{Young diagram $Y(\lambda)$} of a partition $\lambda$ is the set $\{(i,j) : 1\le j \le \lambda_i\}$, drawn as a top-left aligned array of boxes with $\lambda_i$ boxes in row $i$. The \emph{transpose $\lambda'$ of $\lambda$} is defined by $Y(\lambda') = \{(i,j) : (j,i)\in Y(\lambda)\}$. If $\lambda\in L(h,w)$, the \emph{box-complement partition} of $\lambda$ in $L(h,w)$ is defined by $\lambda^\square = (w-\lambda_h, w-\lambda_{h-1}, \ldots, w-\lambda_1)$. 
We write $\mu\subseteq\lambda$ if $Y(\mu)\subseteq Y(\lambda)$ and $\lambda\cap\mu$ for the partition such that $Y(\lambda\cap\mu) = Y(\lambda)\cap Y(\mu)$.
A \emph{skew partition} $\lambda/\mu$ is a pair of partitions such that $\mu\subseteq\lambda$. Its Young diagram is $Y(\lambda/\mu) = Y(\lambda)\setminus Y(\mu)$.
\medskip

Given a tuple of non-negative integers $(\alpha_1, \ldots, \alpha_n)$, let $\mathbf{x}_n^\alpha = x_1^{\alpha_1}\cdots x_n^{\alpha_n}$.
The \emph{monomial symmetric polynomial} $m_\lambda(\mathbf{x}_n)$ indexed by $\lambda$ is the sum $\sum_{\alpha \in S_n.\lambda} \mathbf{x}_n^\alpha$ over all permutations of $(\lambda_1, \ldots, \lambda_{\ell(\lambda)}, 0, \ldots, 0) \in \N_0^n$. The \emph{power sum symmetric polynomial} indexed by $k$ is $p_k(\mathbf{x}_n) = x_1^k + \cdots + x_n^k$. For a partition $\lambda$, set $p_\lambda(\mathbf{x}_n) = \prod_{1\leq i\leq \ell(\lambda)} p_{\lambda_i}(\mathbf{x}_n)$. The \emph{Schur polynomial} indexed by $\lambda$ is 
\[
s_\lambda(\mathbf{x}_n) = \det(x_j^{\lambda_i+n-i})_{1\le i,j\le n}/\det(x_j^{n-i})_{1\le i,j\le n}.
\]
\begin{thm}
    The families $\{m_\lambda(\mathbf{x}_n)  : \lambda\in\Par_{\le n}\}, \{p_\lambda(\mathbf{x}_n)  : \lambda'\in\Par_{\le n}\},$ and \break$\{s_\lambda(\mathbf{x}_n)  : \lambda\in\Par_{\le n}\}$ of symmetric polynomials are $\C$-bases of $\Lambda[\mathbf{x}_{n}]$.
\end{thm}
By taking colimits (or direct limits) of these families, we obtain $\C$-bases \break$\{m_\lambda  : \lambda\in\Par\}$, $\{p_\lambda  : \lambda\in\Par\},$ and $\{s_\lambda  : \lambda\in\Par\}$ of $\Lambda[\mathbf{x}_\infty]$.

We remark that $s_\lambda(1, \stackrel{n}{\ldots}, 1)$ counts the number of \emph{semistandard Young tableaux $T\in\mathrm{SSYT}_n(\lambda)$} in $n$ letters and of shape $\lambda$ --- see \cite[\S7.10]{StanleyEC2} for the necessary definitions. Stanley's \emph{hook-content formula} is
\[
\#\mathrm{SSYT}_n(\lambda) = \prod_{u\in Y(\lambda)}\frac{n+c(u)}{h(u)},
\]
where for $u = (i,j) \in Y(\lambda)$, the number $h(u) = \lambda_i+\lambda'_j-i-j+1$ is its \emph{hook length}, and $c(u) = j-i$ its \emph{content}. The product is taken to be $0$ if $\mathrm{SSYT}_n(\lambda)=\emptyset$.
\medskip

The \emph{Hall scalar product} is the bilinear form $\langle\cdot\,,\cdot\rangle$ on $\Lambda[\mathbf{x}_\infty]$ for which the Schur functions are orthonormal: $\langle s_\lambda, s_\mu \rangle = \delta_{\lambda,\mu}$, where $\delta$ is the Kronecker delta. It follows that every $f \in \Lambda[\mathbf{x}_\infty]$ has a unique expansion $f = \sum \langle f, s_\lambda\rangle s_\lambda$. The basis orthonormal to monomial symmetric functions is denoted $\{h_\lambda  : \lambda\in\Par\}$.

The numbers $c_{\nu\mu}^\lambda = \langle s_\nu s_\mu, s_\lambda\rangle$ are the \emph{Littlewood--Richardson coefficients}. The \emph{skew Schur function} indexed by $\lambda/\mu$ is defined by $s_{\lambda/\mu} = \sum c_{\nu\mu}^\lambda s_\nu$.\medskip

In what follows, we consider alphabets $\mathbf{x}$, $\mathbf{y}$, $\mathbf{z}$ which may be either finite or infinite. 
The \emph{Cauchy identity} states
\[
\sum_{\lambda\in\Par} s_\lambda(\mathbf{x}) s_\lambda(\mathbf{y}) = \prod_{i,j}(1-x_i y_j)^{-1}.
\]
The \emph{addition} and \emph{product} of $\mathbf{x}$ and $\mathbf{y}$ are the alphabets
\[
\mathbf{x}+\mathbf{y} = \{x_1, y_1, x_2, y_2, \ldots\}
\quad\text{and}\quad
\mathbf{xy} = \{x_1y_1, x_1y_2, x_2y_1, \ldots\}.
\]
Given an alphabet $\mathbf{z}$, the map $f\mapsto f[\mathbf{z}]$ is a $\C$-algebra homomorphism \cite[\S2.3]{LoehrRemmel}. The \emph{plethystic addition formula} \cite[\S3.2]{LoehrRemmel} asserts that
\begin{equation}
    s_\lambda[\mathbf{x}+\mathbf{y}] = \sum\nolimits_{\mu} s_\mu(\mathbf{x}) s_{\lambda/\mu}(\mathbf{y}).
\label{eq:plet addition}
\end{equation}
Define the \emph{Kronecker coefficients} $g_{\lambda\mu\nu}$ 
as the constants appearing in
\[
s_\lambda[\mathbf{xy}] =
\sum_{|\mu|=|\nu|=|\lambda|} g_{\lambda\mu\nu}\,s_\nu(\mathbf{x}) s_\mu(\mathbf{y}).
\]
The \emph{Kronecker product} is defined by $(s_\lambda * s_\mu)(\mathbf{x}) = \sum_\nu g_{\lambda\mu\nu}\,s_\nu(\mathbf{x})$ on the Schur basis, and extended bilinearly to $\Lambda[\mathbf{x}]$. Combining these expressions we obtain the \emph{generalised Cauchy identity},
\begin{equation}\label{eq:generalised Cauchy}
s_\lambda[\mathbf{x}\mathbf{y}] = \sum\nolimits_\mu (s_\lambda * s_\mu)(\mathbf{x}) \cdot s_\mu(\mathbf{y}).
\end{equation}

Kronecker coefficients $g_{\lambda\mu\nu}$ are invariant under permutations of $\lambda, \mu, \nu$. A necessary condition for $g_{\lambda\mu\nu}\neq 0$ is $\ell(\nu)\leq\ell(\lambda)\ell(\mu)$. A stronger result states that if $|\lambda|=|\mu|$ then $\max\{\ell(\nu): g_{\lambda\mu\nu}\neq0\}=|\lambda\cap\mu'|$ \cite{Dvir}.

\subsection{\texorpdfstring{$q$}{q}-analogues}
Let $n$ be a non-negative integer. Its $q$-analogue is the $q$-\emph{integer} $[n] = 1+q+\cdots+q^{n-1}$. The $q$-\emph{binomial} is defined as
\[
\qbinom{n}{k} = \frac{[n][n-1]\cdots[n-k+1]}{[k][k-1]\cdots[1]}.
\]

The ring of $\GL_2(\CC)$-characters is 
\[
\langle s_\lambda(x,y) \mid \ell(\lambda)\le2\rangle_{\N_0}.
\]
If $f(x,y)$ is the character of a homogeneous representation, then
we can recover $f(x,y)$ from the specialisation $f(1,q)$. By abuse of notation, we say that $f(1,q)$ is a $\GL_2(\CC)$-character. Note that if $\ell(\lambda)\le2$ then
\[
s_\lambda(1,q) = q^{\lambda_2} s_{(\lambda_1-\lambda_2)} = q^{\lambda_2} [\lambda_1-\lambda_2+1].
\]

All representations in this paper are plethysms of the form $\Sym^m\Sym^n\CC^2$. The character of such a representation is 
\[
s_{(m)}\circ s_{(n)}(1,q) = s_{(m)}(1,q,\ldots,q^n) = \qbinom{n+m}{n},
\]
where $\circ$ denotes plethysm of symmetric polynomials~\cite[Ch.~7, Def.~A2.6]{StanleyEC2}. By the above considerations, there exist non-negative integers 
$a_{m[n]}^k$ such that
\[
\qbinom{n+m}{n} = \sum_{k\ge0} q^{(mn-k+1)/2} \,a_{m[n]}^k \cdot [k].
\]

\subsection{Representations as polynomial rings}\label{sec:preliminaries.poly rings}

Let $\GL_2(\F)$ be the group of $2\times2$ invertible matrices and let $\SL_2(\F)$ be the subgroup of those matrices with determinant $1$.

When $\F=\CC$ they admit a Lie group structure. 
The irreducible highest weight representations of $\SL_2(\C)$ are indexed by $n\in\N$ and classically constructed as the symmetric powers $\Sym^n\CC^2$. 
There are two isomorphic constructions of these symmetric power functors, and distinguishing between them is crucial in the rest of this paper. We write $\Sym_n V$ for the space of invariants
\[
\Sym_n V = (V^{\otimes n})^{S_n}
\]
and $\Sym^n V$ for the space of coinvariants
\begin{equation}\label{eq:covariants}
\Sym^n V = V^{\otimes n}/\langle
\cdots \otimes u\otimes v\otimes \cdots = \cdots \otimes v\otimes u\otimes \cdots
\rangle = \C_n[V],
\end{equation}
where $\C_n[V]$ denotes the space of polynomials of homogeneous degree $n$ whose formal variables are the elements of $V$.
These two constructions are in general non-isomorphic over an arbitrary field $\F$ \cite{McDowellWildon}, and this insight is key in our interpretation of representations as polynomial rings \cite{GMSW, GMSW-fpsac}.\medskip

When $V = \F^2$ we typically choose a basis $\F^2 = \langle x, y\rangle$ and write $\Sym^n \F^2 = \F_n[x,y]$. 
The action of $\GL_2(\F)$ is given by
\[
\left(\begin{smallmatrix}
    \alpha&\beta\\\gamma&\delta
\end{smallmatrix}\right)\!.\,p(x,y)
=p(\alpha x+\gamma y,~\beta x+\delta y)
=(\beta x+\delta y)^n \cdot p\Big(\frac{\alpha x+\gamma y}{\beta x+\delta y},1\Big).
\]
For convenience, we can specialise $y=1$ and work with the space $\F_{\le n}[x]$ of polynomials in one variable with degree at most $n$. The action is then given by
\[
\left(\begin{smallmatrix}
    \alpha&\beta\\\gamma&\delta
\end{smallmatrix}\right)\!.\,p(x)
=(\beta x+\delta )^n \cdot p\Big(\frac{\alpha x+\gamma }{\beta x+\delta }\Big).
\]

We extend this framework to the plethysms $\Sym^m\Sym^n\CC^2$ appearing in Conjecture~\ref{conj:0}, which we choose to interpret as $\Sym_m\Sym^n\F^2$. Combining the constructions and conventions above gives
\[
\Sym_m \Sym^n\F^2 \cong \Sym_m (\F_{\le n}[x]) \cong \Lambda^\F_{\le n}[\mathbf{x}_m],
\]
where $\Lambda^\F_{\le n}[\mathbf{x}_m]$ is the 
$\F$-vector space 
of symmetric polynomials in $\mathbf{x}_m$, with degree at most $n$ in each of $m$ variables.
Following the $\GL_2(\F)$-action throughout the steps of the construction, we arrive at
\begin{equation}\label{matrix action}
\left(\begin{smallmatrix}
    \alpha&\beta\\\gamma&\delta
\end{smallmatrix}\right)
\!.\,f(x_1, \ldots, x_m)
=
\prod_{i=1}^m(\beta x_i+\delta)^n\cdot f\!\left(
    \frac{\alpha x_1+\gamma}{\beta x_1+\delta},
    \ldots,
    \frac{\alpha x_m+\gamma}{\beta x_m+\delta}
    \right).
\end{equation}
Recall $L(h,w) = \{\text{partitions fitting in a $w\times h$ box}\} = \{\lambda\in\Par_{\le h} : \lambda_1 \le w\}$. 
\begin{thm}
    The following are bases of $\Lambda^\F_{\le w}[\mathbf{x}_h]$:
    \begin{enumerate}
        \item $\{m_\lambda(\mathbf{x}_h)  : \lambda \in L(h,w)\}$,
        \item $\{s_\lambda(\mathbf{x}_h)  : \lambda \in L(h,w)\}$.
    \end{enumerate}
\end{thm}
\begin{proof}{\, }
\begin{enumerate}
    \item 
    We know that $\{m_\lambda(\mathbf{x}_h)  : \lambda \in L(h,w)\}$ is a linearly independent set in $\Lambda[\mathbf{x}_h]$. In fact, each $m_\lambda(\mathbf{x}_h)$ is in $\Lambda^\F_{\le w}[\mathbf{x}_h]$ because $\lambda_1 \le w$. Finally, if \break$f = \sum_{\alpha_i \le w~\forall i} c_\alpha \mathbf{x}_h^\alpha \in \Lambda^\F_{\le w}[\mathbf{x}_h]$ then 
    $f = \sum_{\lambda\in L(h,w)} c_\lambda m_\lambda(\mathbf{x}_h),$ showing it is a generating set.
    \item
    We know that $\{s_\lambda(\mathbf{x}_h)  : \lambda \in L(h,w)\}$ is a linearly independent set in $\Lambda^\F[\mathbf{x}_h]$. Since $\{m_\lambda(\mathbf{x}_h)  : \lambda\in\Par_{\le h}\}$ is a basis of $\Lambda^\F[\mathbf{x}_h]$, we can define $K_{\lambda,\mu}$ by $s_{\lambda}(\mathbf{x}_h) = \sum_{\mu} K_{\lambda, \mu}m_\mu(\mathbf{x}_h)$. We have $K_{\lambda,\mu} = 0$ unless $\mu_1 + \cdots + \mu_i \le \lambda_1 + \cdots + \lambda_i$ for all $i = 1, \ldots, h$ and $K_{\lambda,\lambda}=1$ (see \cite[Prop.~7.10.5]{StanleyEC2}). In particular: unless $\mu_1 \le \lambda_1 \le w$. Hence, $\{s_\lambda(\mathbf{x}_h)  : \lambda \in L(h,w)\}$ is in the span of $\{m_\mu(\mathbf{x}_h)  : \mu \in L(h,w)\}$, which is $\Lambda^\F_{\le w}[\mathbf{x}_h]$ by the previous part. Since it is linearly independent and of the same size, it forms a basis.\qedhere
\end{enumerate}
\end{proof}

\section{Consequences of Theorem~\ref{thm:main} for Gaussian coefficients}\label{sec:Troyka}

Before we move on to studying the two maps which realise the generalised Foulkes conjecture for $\SL_2(\C)$, we give an overview of some corollaries of Theorem~\ref{thm:main} related to Gaussian coefficients.

In an attempt to tackle the generalised Foulkes conjecture, Bergeron~\cite{BergeronSlides} posed a problem about the coefficients of $\qbinom{c+d}{c}-\qbinom{a+b}{a}$ as a polynomial in $q$ being non-negative. Theorem~\ref{thm:main} gives a proof under the divisibility condition.

\setcounter{section}{1}
\setcounter{de}{4}
\begin{cor}
Let $a, b, c, d \in \N$, such that $ab=cd$, and assume $a=\min(a,b,c,d)$. 
If $a$ divides either $c$ or $d$, then $\qbinom{c+d}{c} - \qbinom{a+b}{a} \in \N_0[q]$.
\end{cor}

In \cite{Zanello}, Zanello further conjectured that this polynomial is unimodal. Again, Theorem~\ref{thm:main} confirms this conjecture under divisibility conditions.

\begin{cor}
Let $a, b, c, d \in \N$, such that $ab=cd$ and $a=\min(a,b,c,d)$. If $a$ divides either $c$ or $d$, then $\qbinom{c+d}{c} - \qbinom{a+b}{a}$ is a unimodal polynomial.
\end{cor}

\begin{proof}[Proof of Corollaries \ref{cor: Bergeron} and \ref{cor: Zanello}]
		Consider the injection of $\SL_2(\C)$-representations $$\Sym^b\Sym^a\CC^2 \hookrightarrow \Sym^d\Sym^c\CC^2.$$
		Since $\SL_2(\C)$ is semisimple, there is a complement $C$ satisfying
		\[
		\Sym^b\Sym^a\CC^2 \oplus C \cong \Sym^d\Sym^c\CC^2.
		\]
		The degrees of both sides match, that is $ab=cd$, so $C$ lifts to a complement of $\GL_2(\C)$-representations. 
		The $\GL_2(\C)$-character of $C$ is $\qbinom{c+d}{c} - \qbinom{a+b}{a}$, and thus the coefficients are positive. This proves Corollary \ref{cor: Bergeron}.
	    To finish the proof of Corollary \ref{cor: Zanello}, recall that the coefficients of a polynomial $f\in\N_0[q]$ form a symmetric and unimodal sequence if and only if $f$ is the character of a $\GL_2(\C)$-representation of homogeneous degree \cite[Thm.~2.1]{StanleySL2}.
\end{proof}

Interpreting the $q$-binomial coefficients as generating functions for partitions fitting inside a rectangle, it is clear that for any $a,b,c,d$, a necessary condition for $q$-positivity of $\qbinom{c+d}{c} - \qbinom{a+b}{a}$ is $ab\le cd$. Using the \emph{$q$-Pascal identity}
\[
\qbinom{n}{k} - \qbinom{n-1}{k} = q^{n-k}\qbinom{n-1}{k-1},
\]
Corollary~\ref{cor: Bergeron} immediately generalises to the case where the equality $ab=cd$ is relaxed to the condition $ab\le cd$. Note that this is a necessary condition, since $ab$ is the degree of the polynomial $\qbinom{a+b}{a}$.

\begin{cor} \label{cor: Troyka}
Let $a,b,c,d\in\N$, such that $ab\le cd$ and $a=\min(a,b,c,d)$. If $a$ divides either $c$ or $d$ then $\qbinom{c+d}{c}-\qbinom{a+b}{a}\in\N_0[q]$.
\end{cor}
\begin{proof}
    By Corollary~\ref{cor: Bergeron} for parameters $a,\frac{cd}{a},c,d$,
    \[
    \qbinom{c+d}{c}-\qbinom{a+\frac{cd}{a}}{a}\in\N_0[q].
    \]
    Now applying the $q$-Pascal identity $\frac{cd}{a}-b\ge0$ times, we obtain
    \[
    \qbinom{a+\frac{cd}{a}}{a}-\qbinom{a+b}{a}\in\N_0[q].
    \]
    Adding the two $q$-positive expressions above, we deduce $\qbinom{c+d}{c} - \qbinom{a+b}{a}\in\N_0[q]$.
\end{proof}

In \cite{Troyka}, Troyka conjectures that this more general result holds without the divisibility condition.

\begin{conjecture}[Troyka]\label{conj:Troyka}
Let $a,b,c,d\in\N$, such that $ab\le cd$ and $a=\min(a,b,c,d)$. Then $\qbinom{c+d}{c}-\qbinom{a+b}{a}\in\mathbb{N}_0[q]$.
\end{conjecture}

Troyka further notes that for certain choices of parameters $a,b,c,d$, this conjecture can be deduced from the $\SL_2(\C)$ case of Conjecture~\ref{conj:A} using Butler's result on $q$-unimodality \cite{Butler}. 
\setcounter{section}{3}
\setcounter{de}{0}
\begin{p}[Butler]\label{p:Butler}
    For all $n\in\N$ the sequence $\qbinom{n}{0}, \qbinom{n}{1}, \ldots, \qbinom{n}{n}$ is \break$q$-unimodal. That is,
\begin{equation}\label{eq:Butler}
\qbinom{n}{k} - \qbinom{n}{k-1} \in \N_0[q] \quad\text{for all }0\le k \le \Big\lfloor \frac{n}{2}\Big\rfloor.
\end{equation}
\end{p}
However, since the generalised Foulkes conjecture remained open at the time, these deductions were necessarily conditional. With Conjecture~\ref{conj:A} now settled under the divisibility conditions by Theorem~\ref{thm:main}, we can make these deductions rigorous and unconditional, further generalising Corollary~\ref{cor: Troyka}.

\begin{p}\label{prop:BestTroyka}
    Let $a,b,c,d$ be positive integers, such that $a$ is the smallest. If there exists a pair of integers $x\le y$ such that
    \begin{enumerate}
        \item $a\le x\le \min(c,d)$,
        \item $x+y\le c+d$,
        \item $xy\ge ab$,
        \item $a$ divides either $x$ or $y$,
    \end{enumerate}
    then $\qbinom{c+d}{c}-\qbinom{a+b}{a}\in\N_0[q]$.
\end{p}

\begin{proof}
    Let $x,y$ be such a pair, and without loss of generality assume $c\le d$, so $c=\min(c,d)$. Since $1\leq a\le x\le c\le \lfloor\frac{c+d}{2}\rfloor$, Proposition~\ref{p:Butler} gives
    \[
    \qbinom{c+d}{c}-\qbinom{c+d}{x}=\sum_{k=1}^{c-x}\qbinom{c+d}{x+k}-\qbinom{c+d}{x+k-1}\in\N_0[q].
    \]
    By condition (2) and the $q$-Pascal identity, we also have
\begin{multline*}
    \qbinom{c+d}{x}-\qbinom{x+y}{x}=\sum_{k=1}^{c+d-x-y}\qbinom{x+y+k}{x}-\qbinom{x+y+k-1}{x}\\=\sum_{k=1}^{c+d-x-y}q^{y-k}\qbinom{x+y+k-1}{x-1}\in\N_0[q].
\end{multline*}
    Conditions (3) and (4) allow us to apply Corollary~\ref{cor: Troyka} for the parameters $(a,b,x,y)$, so that
        \[
    \qbinom{x+y}{x}-\qbinom{a+b}{a}\in\N_0[q],
    \]
    as well. Adding the three $q$-positive expressions above finishes the proof.
\end{proof}

\section{The \texorpdfstring{$k$}{k}-fold plethystic substitution}
\label{sec:main}
	In this section we recall the construction of the $k$-fold map from \cite{GMSW-fpsac}, and show that it realises a $\GL_2(\F)$-invariant map 
\[
    \Sym_{kd}\Sym^a\F^2\to\Sym_{d}\Sym^{ka}\F^2.
\]
    When $\F=\C$, letting $kd=b$ and $ka=c$ turns the $k$-fold map into a candidate for Conjecture~\ref{conj:A} for $\SL_2(\C)$.
    
    We start by defining the product of an alphabet $\mathbf{x}$ with a scalar $k\in\N$ as the alphabet
    \[
    k\mathbf{x} = \{x_1, \stackrel{k}{\ldots} , x_1, x_2, \stackrel{k}{\ldots} , x_2, \ldots\}.
    \]
	\begin{de}\label{de:k-fold map}
		The \emph{$k$-fold plethystic substitution}, or the \emph{$k$-fold map}, is defined as
		\begin{align*}
			\kappa^k_{\infty,d} : \Lambda^{\F}[\mathbf{x}_{kd}] &\to  \Lambda^{\F}[\mathbf{y}_{d}] \\
			f(\mathbf{x}_{kd}) &\mapsto
			f[k\mathbf{y}_d] = f(y_1, \stackrel{k}{\ldots}, y_1,\ldots,y_d,\overset{k}{\ldots}\,,y_d).
		\end{align*}
	\end{de}
	Equivalently, the $k$-fold map 
	sends each $x_i$ to $y_{\lceil i/k\rceil}$ for all $i$. It is an $\F$-algebra homomorphism. The definition extends to infinite alphabets.
    \begin{p}\label{p:well-defnd}
        The image of $\Lambda^\F_{\leq a}[\mathbf{x}_{kd}]$ under the $k$-fold map lies inside $\Lambda^\F_{\leq ka}[\mathbf{y}_d]$.	
    \end{p}
    \begin{proof}
        Since the $k$-fold map is linear, it suffices to consider the image of a single monomial:
		\begin{multline*}
			\Big(x_1^{a_1} x_2^{a_2} \cdots x_{k}^{a_k}\Big)
			\Big(x_{k+1}^{a_{k+1}} \cdots x_{2k}^{a_{2k}}\Big)
			\cdots 
			\Big(x_{kd-k+1}^{a_{kd-k+1}} \cdots x_{kd}^{a_{kd}}\Big)
			\mapsto \prod_{i=1}^d y_i^{a_{(i-1)k+1} + \cdots + a_{ik}}.
		\end{multline*}
		Since $a_j\leq a$ for all $j$, the degree of $y_i$ in the image is
		\[
		{a_{(i-1)k+1} + \cdots + a_{ik}} \le k \cdot a,
		\]
		for each $i = 1, \ldots, d$, as claimed.
    \end{proof}
    We denote by $\kappa_{a,d}^k$ the restriction $\Lambda^\F_{\leq a}[\mathbf{x}_{kd}] \rightarrow  \Lambda^\F_{\leq ka}[\mathbf{y}_{d}]$ of the $k$-fold map. By abuse of language, we refer to $\kappa_{a,d}^k$ as the $k$-fold map.
	\begin{p}
    The $k$-fold map $\kappa_{a,d}^k$ is $\GL_2(\F)$-equivariant.
	\end{p}
	\begin{proof}
		We compute directly that the $k$-fold map commutes with the action of $\GL_2(\F)$. Apply first the $k$-fold map and then the action of a matrix as in \eqref{matrix action} to get
        \[
    \prod_{i=1}^d\left((\beta y_i+\delta)^a\right)^k\cdot f\!\left(
    \frac{\alpha y_1+\gamma}{\beta y_1+\delta},
    \stackrel{k}{\ldots}\,,
    \frac{\alpha y_1+\gamma}{\beta y_1+\delta}, \ldots, \frac{\alpha y_d+\gamma}{\beta y_d+\delta},
    \stackrel{k}{\ldots}\,,
    \frac{\alpha y_d+\gamma}{\beta y_d+\delta}
    \right)\,.
        \]
        Now apply first the matrix action and then the $k$-fold map, which clearly gives the same result.
	\end{proof}

    It remains to show that the $k$-fold map $\kappa_{a,d}^k$ is injective for $\F=\C$ and $a\le d$. This is a corollary of Theorem~\ref{thm:main} after we have established Theorem~\ref{duality}.

    \subsection*{The \texorpdfstring{$k$}{k}-fold map in literature}
    We remark that while the consideration of the $k$-fold map as a description of $\SL_2(\CC)$-injection of plethysms is novel, the map itself has previously appeared in the literature in two independent contexts we outline below.

    From the combinatorial point of view, it is more natural to consider the map $\kappa_{\infty,\infty}^k$ on symmetric functions, where the number of variables and their degrees are unrestricted. The map admits a simple description on the power sum basis:
    \[
    \kappa_{\infty,\infty}^k:p_n\mapsto kp_n,
    \]
    and so is immediately seen injective. Moreover, its structure constants in the Schur basis and in the monomial basis exhibit several beautiful properties of combinatorial nature which have been studied by \cite{Baker, Brenti}, some of which we review in \S\ref{sec:matrix}.

    On the other hand, the $k$-fold map $\kappa_{a,d}^k$ with finite parameters appeared in the condensed matter physics literature. In \cite{GEA}, the authors reduce the question of whether the rank of the particle entanglement spectrum of Laughlin states matches the number of quasi-hole excitations to showing that the $k$-fold map is injective. This is now settled in Corollary~\ref{cor:physics} of this article.

\section{The Foulkes--Howe map and the \texorpdfstring{$k$}{k}-fold map}\label{sec: proof of duality}

Fix $k$ and $a$. This section is dedicated to the proof of Theorem~\ref{duality}. We show that the Foulkes--Howe map $\alpha_a$
\[
 \Sym^{kd}\Sym^aU \leftarrow
\Sym^d\Sym^{ka}U
\]
as constructed in \cite[(1.3b)]{RSWY} is the dual of the $k$-fold map $\kappa_{a,d}^k$. Hence for $d\ge a$ we deduce $\alpha_d$ is surjective if and only if $\kappa_{a,d}^k$ is injective. Either of these maps realises the sought proof of Conjecture~\ref{conj:A} for $\SL_2(\C)$ (under divisibility conditions) via Theorem~\ref{thm:main}.\medskip

One of the main goals of \cite{RSWY} is to understand the minimal generators of the defining ideal 
of the varieties $X_{(k^a)}$ defined below, as well as their higher syzygy modules. These varieties (known as coincidence root loci \cite{AC}) are
defined as
$$X_{(k^a)}=\{[F]\in\mathbb{P}^m\mid F=L_1^{k}L_2^{k}\dots L_a^{k},\textrm{ where each }L_i\textrm{ is a linear form}\}.$$
Part of the main result is the statement that $I(X_{(k^a)})$ contains no equations of degree $\leq a$, and this relies on the interpretation of $X_{(k^a)}$ as a linear projection of a Veronese variety
\begin{center}
    \begin{tikzcd}[row sep=2.5em]
 & \mathbb{P}^N
 \arrow[dr, dashed] \\
\mathbb{P}^a \arrow[ur,hook] \arrow{rr}{\nu_k} && \mathbb{P}^{ka}
\end{tikzcd}
\end{center}
where $N=\binom{k+a}{k}-1$. By interpreting the top space as the projectivisation of a degree component of the coordinate ring of $\mathbb{P}^a$, one can induce the left map. The bottom map is an exponentiation map defined below. The dashed map is induced by the inclusion $\Sym^{ka}\C^2\hookrightarrow\Sym^k\Sym^a\C^2$.
See \cite{AC14} for more details on this construction. 
Given that this approach is so vastly different from ours, we require a step-by-step translation of the set-up.\medskip

Even though the authors of~\cite{RSWY} work only over $\C$, we chose to translate everything to a field-independent setting. Although a priori this appears to complicate the argument, it in fact allows us to keep track of the correct dualities implicit in each step of the construction. Moreover, we will sometimes distinguish between a vector space $V$ and the canonically isomorphic $V^{\star\star}$, to treat its elements as functions mapping elements of $V^\star$ to $\F$.\medskip

Let $U = \F^2 = \langle x,y\rangle_\F$ be the natural representation of $\SL_2(\F)$, and denote the dual $U^\star = \langle X, Y\rangle_\F$.
Here, $X(\alpha x+\beta y) = \alpha$ and $Y(\alpha x+\beta y) = \beta$. 
The main objects of study are binary forms of degree $a$, which are polynomials of the form
\[
F = f_0Y^a + f_1XY^{a-1} + \cdots + f_a X^a.
\]
That is, $F \in \F_a[U^\star] = \Sym^aU^\star$. These forms are studied up to scaling, and so they form a projective space which the authors of \cite{RSWY} denote by $\mathbf{P}^a$. Write $[F]\in\mathbf{P}^a$ for the element in the projective space corresponding to $F$.   Following our convention from \S\ref{sec:preliminaries.poly rings}, we let $Y=1$ and write forms as polynomials in $\F_{\le a}[X]$ instead.

Next, they consider the map
\begin{align*}
    \nu_k : \mathbf{P}^a &\to \mathbf{P}^{ka}\\
    [F] &\mapsto [F^k],
\end{align*}
where 
\(
F^k = (f_0 + f_1X + \cdots + f_aX^a)^k
\)
is the usual exponentiation in $\F[X]$. 
Note that the coefficient of $X^i$ in $F^k$ is
\begin{equation}
\mathrm{coeff}_{X^i}(F^k) = \sum_{\substack{\beta_1+\cdots+\beta_k=i\\\beta_j\le a ~\forall j}} f_{\beta_1}f_{\beta_2}\cdots f_{\beta_k}.
\label{eq:coeff F^k}
\end{equation}

Dual to this map, there is a map of coordinate rings
\[
\mathcal{O}_{\mathbf{P}^a} \xleftarrow{\nu_k^\star} \mathcal{O}_{\mathbf{P}^{ka}}
\]
which is given by pre-composing with $\nu_k$. The coordinate ring of a vector space seen as a variety is simply $\mathcal{O}_V = \F[V^\star] = \Sym^\bullet V^\star$, where the symbol $\Sym$ means that the tensor factors commute (cf.~\eqref{eq:covariants}). Since $\mathbf{P}^a$ is the projectivisation of $\Sym^aU^\star$ and since $(\Sym^aU^\star)^\star = \Sym_aU^{\star\star}$ \cite{McDowellWildon}, we have defined (up to scaling) a map
\[
\Sym^\bullet \Sym_aU^{\star\star} \leftarrow
\Sym^\bullet \Sym_{ka}U^{\star\star}.
\]
It is important now to see how this map behaves on a basis. A basis of $\Sym_aU^{\star\star}$ is given by the elements
\[
\sum_{\sigma\in S_a}
\sigma( X^\star\otimes\stackrel{i}{\cdots}\otimes X^\star 
\otimes
Y^\star\otimes\stackrel{a-i}{\cdots}\otimes Y^\star )
= \mathrm{coeff}_{X^i}(-)
\]
as $i$ ranges over $0, \ldots, a$. Applying this element to a form $F\in\F_{\le a}[X]$ extracts the coefficient $f_i$ of $X^i$.
Introducing a different alphabet $X_1, Y_1, X_2, Y_2, \ldots$ for each copy of $\Sym_aU^{\star\star}$ in $\Sym^\bullet \Sym_aU^{\star\star}$, 
Equation~\eqref{eq:coeff F^k} translates to 
\begin{align*}
\nu_k^\star :\mathrm{coeff}_{X^i} &\mapsto \sum_{\substack{\beta_1+\cdots+\beta_k=i\\\beta_j\le a ~\forall j}} \mathrm{coeff}_{X_1^{\beta_1}}\mathrm{coeff}_{X_2^{\beta_2}}\cdots \mathrm{coeff}_{X_k^{\beta_k}}
\\&= 
\sum_{\lambda\in L(k,a)\cap\Par(i)} \binom{k}{m_0(\lambda),m_1(\lambda),\ldots,m_a(\lambda)} \cdot \mathrm{coeff}_{\mathbf{X}^{\lambda}}\,,
\end{align*}
where $\mathbf{X}^{\lambda} = X_1^{\lambda_1}\cdots X_k^{\lambda_k}$. The multiplication in $\Sym^\bullet \Sym_aU^{\star\star}$ is commutative and the domain of the function on the right-hand side is a space of symmetric polynomials, and hence $\mathrm{coeff}_{\mathbf{X}^{\lambda}}$ is well defined.\medskip

In a subsequent step, the authors translate this map to a map between sheaf cohomology groups. To do this, they invoke a classical result on global sections of projective spaces \cite[\S14.1.2]{Vakil}: given a vector space $V$, let $\mathbb{P}(V)$ be the projective space parametrising $1$-dimensional
quotients of $V$; then its space of twisted global sections is
\[
H^0(\mathbb{P}(V), \mathcal{O}_{\mathbb{P}(V)}(d)) = \Sym^dV.
\]
In our case, $\mathbf{P}^a$ is the projective space parametrising subspaces of $\Sym^aU^\star$, and hence $\mathbf{P}^a = \mathbb{P}(\Sym_aU^{\star\star})$ as parametrising quotients. 
Therefore the map $\nu_k^\star$ induces a map as follows:
\[
\begin{tikzcd}[row sep=tiny]
\Sym^{kd}\Sym_{a}U^{\star\star}                                                        &  & \Sym^d\Sym_{ka}U^{\star\star} \arrow[ll, "\alpha_d"']                                                           \\
{H^0(\mathbf{P}^{a}, \mathcal{O}_{\mathbf{P}^{a}}(kd))} \arrow[u, Rightarrow, no head] &  & {H^0(\mathbf{P}^{ka}, \mathcal{O}_{\mathbf{P}^{ka}}(d))} \arrow[ll, "\alpha_d"] \arrow[u, Rightarrow, no head]
\end{tikzcd}
\]
A basis of the top-left space is given by
\(\{\mathrm{coeff}_{\mathbf{X}^\lambda}\mid\lambda\in L(kd,a)\}\), and similarly \(\{\mathrm{coeff}_{\mathbf{X}^\mu}\mid\mu\in L(d,ka)\}\) for the top-right space. 
The map $\alpha_d$ is
\begin{equation}
    \label{eq:alpha in coordinates}
\alpha_d : \mathrm{coeff}_{\mathbf{X}^\mu} \mapsto
\prod_{i=1}^{\ell(\mu)}\sum_{\lambda\in L(k,a)\cap\Par(\mu_i)} \binom{k}{m_0(\lambda),m_1(\lambda),\ldots,m_a(\lambda)} \cdot\mathrm{coeff}_{\mathbf{X}^\lambda}\,.
\end{equation}

On the other hand, using the realisation of plethysm of symmetric powers from \cite{GMSW}, the $k$-fold map fits in the following diagram:
\[
\begin{tikzcd}[row sep=tiny]
\Sym_{kd}\Sym^{a}U^{\star} \arrow[rr, "\kappa_{a,d}^k"]                                        &  & \Sym_d\Sym^{ka}U^{\star}                                          \\
{\Lambda_{\le a}[\mathbf{X}_{kd}]} \arrow[rr, "\kappa_{a,d}^k"'] \arrow[u, Rightarrow, no head] &  & {\Lambda_{\le ka}[\mathbf{X}_{d}]} \arrow[u, Rightarrow, no head]
\end{tikzcd}
\]
Its dual is the map
\begin{align*}
(\kappa_{a,d}^k)^\star : \mathrm{coeff}_{\mathbf{X}^\mu} \mapsto& 
\Big(
\Lambda_{\le a}[\mathbf{X}_{kd}] \ni f
\mapsto \mathrm{coeff}_{\mathbf{X}^\mu} \big(f[k\mathbf{X}_d]\big) \in \F
\Big)\\
=& 
\prod_{i=1}^{\ell(\mu)}\sum_{\substack{\beta^{(i)}_1+\cdots+\beta^{(i)}_k = \mu_i\\ \beta_j^{(i)}\le a~\forall j}} \mathrm{coeff}_{\mathbf{X}^{\beta^{(i)}}}\,,
\end{align*}
which recovers~\eqref{eq:alpha in coordinates} by regrouping terms as above. (This can also be seen by dualising the formula in Proposition~\ref{p:m to m}.) This shows that the Foulkes--Howe map is the dual of the $k$-fold map, completing the proof of Theorem~\ref{duality}. A summary of the maps and dualities used in this section can be found in the following diagram.
\[\begin{tikzcd}[row sep=small, column sep = small]
	&& \F & \\
	& {\Sym^aU^\star} & {\Sym^{ka}U^\star} \\
	{\Sym^\bullet\Sym_aU^{\star\star}} &&& {\Sym^\bullet\Sym_{ka}U^{\star\star}} \\
	{\Sym^{kd}\Sym_aU^{\star\star}} &&& {\Sym^d\Sym_{ka}U^{\star\star}} \\
	{\Sym_{kd}\Sym^aU^\star} &&& {\Sym_d\Sym^{ka}U^\star} \\
	& {\Sym^{kd}\Sym_aU^{\star\star}} & {\Sym^d\Sym_{ka}U^{\star\star}} \\
	& \F
	\arrow["{{{\nu_k^\star}\mathrm{coeff}_{X^i}}}", from=2-2, to=1-3]
	\arrow["\nu_k"', from=2-2, to=2-3]
	\arrow["\mathrm{coeff}_{X^i}"', from=2-3, to=1-3]
	\arrow[""{name=0, anchor=center, inner sep=0}, "{\nu_k^\star}"', from=3-4, to=3-1]
	\arrow[""{name=1, anchor=center, inner sep=0}, "{{{\alpha_d}}}"{description}, from=4-4, to=4-1]
	\arrow[""{name=2, anchor=center, inner sep=0}, "{{\kappa_{a,d}^k}}"', from=5-1, to=5-4]
	\arrow["\mathrm{coeff}_{\mathbf{X}^\mu}"', from=6-2, to=7-2]
	\arrow["{{\alpha_d}}"', from=6-3, to=6-2]
	\arrow["{{{\mathrm{coeff}_{\mathbf{X}^\mu}\circ\kappa_{a,d}^k}}}", from=6-3, to=7-2]
	\arrow[shorten >=.5em, shorten <=.5em, Leftarrow, from=1, to=0]
	\arrow["\star", shorten >=.5em, shorten <=.5em,  Leftrightarrow, from=1, to=2]
\end{tikzcd}\]

To conclude the proof of Conjecture~\ref{conj:A} for $\SL_2(\C)$, we invoke Theorem~\ref{thm:main}.
\setcounter{section}{1}
\begin{thm}[{\cite[Thm.~1.4]{RSWY}}]
    For $\F=\C$, the Foulkes--Howe map $\alpha_d$ is surjective for $d\ge a$ and injective for $d \le a$.
\end{thm}
\setcounter{section}{5}
\setcounter{de}{0}
A key step in the proof of the above theorem is to show a $\GL_2(\C)$ isomorphism between two complexes \cite[Theorem 3.1]{RSWY}. We remark that this is also a consequence of \cite[Theorem 1.2]{GMSW} over an arbitrary field $\F$ and in a slightly more general setting. Another key step is to show that certain cohomology groups of these complexes vanish.\medskip

This settles \cite[Conjecture 7.3]{GMSW-fpsac}.
\setcounter{section}{1}
\setcounter{de}{2}
\begin{cor}\label{cor:k-fold-injective}
  If $\F=\CC$ then the $k$-fold map
  $\kappa_{a,d}^k$ is injective for $a\le d$ and surjective for $a \ge d$.
\end{cor} 
\begin{prob}\label{prob:combinatorial proof}
    Establish Corollary~\ref{cor:k-fold-injective} directly and combinatorially.
\end{prob}
\setcounter{section}{5}
\setcounter{de}{0}

Corollary~\ref{cor:k-fold-injective} settles a conjecture in condensed matter physics; see~\cite[\S{IV}]{GEA}. We point the reader to the original article for the necessary context.
\begin{cor}\label{cor:physics}
    The particle entanglement spectrum of the general Laughlin state is rank saturated.
\end{cor}

\section{The structure constants of the \texorpdfstring{$k$}{k}-fold map}\label{sec:matrix}
Fix $k\in\mathbb{N}$. In this section we review and give new results on the (infinite) matrix of the $k$-fold map
$\kappa^k_{\infty,\infty} : \Lambda^\C[\mathbf{x}_\infty] \ni f(\mathbf{x}_\infty) \mapsto f[k\mathbf{y}_\infty]$
with respect to selected bases of the $\C$-algebra of symmetric functions. 
The hope of the authors is that some of the properties highlighted in this section can one day be used to solve Open Problem~\ref{prob:combinatorial proof} of establishing the injectivity of $\kappa_{a,d}^k$ directly.
The matrices exhibit rich structure and properties that are of independent interest; for numerical examples see
Examples~\ref{eg:D(2)} and~\ref{eg:B(2)}.

\subsection{The power sum basis matrix}
In the power sum basis, we have $\kappa_{\infty,\infty}^k p_\lambda = k^{\ell(\lambda)} p_\lambda$. That is, the matrix $A(k)$ of the structure coefficients is an infinite diagonal matrix in which the entry indexed by $\lambda$ is $k^{\ell(\lambda)}$.

\begin{note}
    This matrix is not useful for the purposes of establishing injectivity of the $k$-fold map, since there is no basis of $\Lambda_{\le n}[\mathbf{x}_m]$ in which all elements are of the form $p_\lambda(\mathbf{x}_m)$, in general.
\end{note}

\subsection{The Schur basis matrix}
The Schur basis matrix of the $k$-fold map on symmetric functions was studied in Baker's PhD thesis~\cite[\S3]{Baker}. 
\begin{de}
    For each $\lambda, \mu \in \Par$, define
    $d_\lambda^\mu(k)$ via $s_\lambda[k \mathbf{y}_\infty] = \sum_\mu d_\lambda^\mu(k) s_\mu(\mathbf{y}_\infty)$.
\end{de}

We study $D(k) = \big( d_\lambda^\mu(k)\big)_{\lambda,\mu\in\Par}$ assuming the lexicographic order on $\Par$ throughout.
Observe first that the matrix is block-diagonal.
\begin{lem}
    If $|\lambda|\ne|\mu|$ then $d_\lambda^\mu(k) = 0$.
\end{lem}
\begin{proof}
    The $k$-fold map preserves the total degree of monomials.
\end{proof}
In fact, since $p_\lambda[k\mathbf{y}_\infty] = k^{\ell(\lambda)}p_\lambda(\mathbf{y}_\infty)$, the matrix $D(k)$ is diagonalised by the matrix $X = (\chi^\lambda(\mu))_{\lambda,\mu\in\Par}$ of character values of the symmetric groups. More precisely, $X^{-1}\cdot D(k)\cdot X = \mathrm{diag}(1,k,k,k^2,\ldots) = \mathrm{diag}(k^{\ell(\lambda)})_{\lambda\in\Par}$.\medskip

We give several formulas for $d_\lambda^\mu(k)$. 
The first formula expresses $d_\lambda^\mu(k)$ in terms of $k$-multi-Littlewood--Richardson coefficients, defined by
$c_{\bm{\nu}}^\lambda = \langle s_{\nu^{(1)}}\cdots s_{\nu^{(k)}}, s_\lambda\rangle$, where $\bm{\nu}=(\nu^{(1)},\ldots,\nu^{(k)})\in\Par^k$. It is partially discussed in \cite[\S3.2.1]{Baker}.
\begin{p}\label{p:coefs}
    For all $\lambda, \mu\in\Par$, we have \(
    d_\lambda^\mu(k) = \sum_{\bm{\nu} \in \Par^k} c_{\bm{\nu}}^\lambda c_{\bm{\nu}}^\mu
    \).
\end{p}
\begin{proof}
    By the plethystic addition formula \eqref{eq:plet addition},
    \begin{align*}
        \textstyle
        s_\lambda\big[k \mathbf{y}_\infty\big] =
        s_\lambda\big[\mathbf{y}_\infty+\stackrel{k}{\cdots}+\mathbf{y}_\infty\big]  &= \sum_{\theta^{(1)}\subseteq\theta^{(2)}\subseteq\cdots\subseteq\theta^{(k)}=\lambda} s_{\theta^{(1)}}s_{\theta^{(2)}/\theta^{(1)}}\cdots s_{\lambda/\theta^{(k-1)}} (\mathbf{y}_\infty). 
    \end{align*}
    By definition, $s_{\theta^{(i+1)}/\theta^{(i)}} = \sum_{\nu} c^{\theta^{(i+1)}}_{\theta^{(i)}\nu} s_\nu$, and hence
    \begin{align*}
        \textstyle
        s_\lambda\big[k \mathbf{y}_\infty\big]  &= \sum_{\substack{\theta^{(1)}\subseteq\theta^{(2)}\subseteq\cdots\subseteq\theta^{(k)}=\lambda\\\nu^{(2)},\ldots,\nu^{(k)}}} 
        c_{\theta^{(1)}\nu^{(2)}}^{\theta^{(2)}}
        c_{\theta^{(2)}\nu^{(3)}}^{\theta^{(3)}} \cdots 
        c_{\theta^{(k-1)}\nu^{(k)}}^{\lambda}
        s_{\nu^{(1)}}s_{\nu^{(2)}}\cdots s_{\nu^{(k)}} (\mathbf{y}_\infty). 
    \end{align*}
    Using associativity and renaming $\theta^{(1)}$ to $\nu^{(1)}$, we get
    \begin{align*}
        \textstyle
        s_\lambda\big[k \mathbf{y}_\infty\big]  &= \sum_{\nu^{(1)},\nu^{(2)},\ldots,\nu^{(k)}} 
        c_{\bm{\nu}}^{\lambda}
        s_{\nu^{(1)}}s_{\nu^{(2)}}\cdots s_{\nu^{(k)}} (\mathbf{y}_\infty).
    \end{align*}
    Comparing the coefficients at $s_\mu$ on both sides gives the desired expression.
\end{proof}

To interpret $d_\lambda^\mu(k)$ combinatorially (Proposition~\ref{thm:Schur count}), we need to borrow a few definitions from~\cite[\S3.1]{SXP}. A \emph{$k$-multi-tableau} $\mathbf{T}=(T^{(1)},\ldots,T^{(k)})$ is a $k$-tuple of semistandard Young tableaux. The \emph{$k$-multi-shape} of a $k$-multi-tableau $\mathbf{T}$ is the element $\mathbf{\nu}\in\Par^k$ such that $\nu^{(i)}$ is the shape of $T^{(i)}$. The \emph{content} of $\mathbf{T}$ is the sum of the contents of the $T^{(i)}$, where the sum of compositions is taken entry-wise. The \emph{word} of $\mathbf{T}$ is the concatenation $w(T^{(1)})\dots w(T^{(k)})$, where given a semistandard Young tableau $T$ we write $w(T)$ for the (left-to-right, top-to-bottom) row reading word. We say $\mathbf{T}$ is \emph{lattice} if $w(\mathbf{T})$ is a lattice word.

By \cite[Definition 3.6 and Lemma 6.1]{SXP}, the $k$-multi-Littlewood--Richardson coefficient $c_{\bm\nu}^\lambda$ counts the number of lattice $k$-multi-tableaux of multi-shape $\bm{\nu}$ and content $\lambda$. Together with Proposition~\ref{p:coefs} we obtain a proof of Proposition~\ref{thm:Schur count}.\medskip

Next, we explore the symmetries of $D(k)$.
\begin{lem}\label{lem:symmetries}
    The coefficients $d_\lambda^\mu(k)$ exhibit the following symmetries:
    \begin{enumerate}
        \item $d_\lambda^\mu(k) = d_\mu^\lambda(k)$,
        \item $d_\lambda^\mu(k) = d_{\lambda'}^{\mu'}(k)$,
        \item $d_\lambda^\mu(k) = d_{\lambda^\square}^{\mu^\square}(k)$, where the box-complement of $\lambda$ (resp.~$\mu$) is taken inside $L(b,a)$ (resp.~$L(d,c)$), for all choices of $a, b, c, d$ such that $c/a = b/d = k$.
    \end{enumerate}
\end{lem}
\begin{proof}~
    \begin{enumerate}
        \item This follows by the symmetry of $\lambda$ and $\mu$ in Proposition~\ref{p:coefs}.
        \item This follows from the equality $c^{\lambda}_{\mu,\nu} = c^{\lambda'}_{\mu',\nu'}$ of Littlewood--Richardson coefficients (obtained by taking the $\omega$ involution on the expansion of $s_\mu s_\nu$), and Proposition~\ref{p:coefs}. 
        \item This follows from the duality $c^{\lambda}_{\mu,\nu} = c^{\lambda^\square}_{\mu^\square,\nu^\square}$ of Littlewood--Richardson coefficients (see \cite{BriandRosas}) and Proposition~\ref{p:coefs}. \qedhere
    \end{enumerate}
\end{proof}
Note that part (1) states that the matrix $D(k)$ is symmetric, and (2) says that the blocks of $D(k)$ are invariant under transpositions across their anti-diagonal. The symmetry given by (3) is more subtle. 
Together, these symmetries endow $D(k)$ with a rich combinatorial structure.\medskip

The Cauchy identity gives another formula for $d_\lambda^\mu(k)$. This result is similar in spirit to the discussion in~\cite[\S3.2.1]{Baker}.
\begin{p}
    For all $\lambda, \mu\in\Par$, we have
    $$d_\lambda^\mu(k) = \Big\langle
    s_\lambda(\mathbf{y}_\infty)s_\mu(\mathbf{z}_\infty), ~  \prod_{i,j} (1-y_jz_i)^{-k} \Big\rangle.$$
\end{p}
\begin{proof}
    Consider the following computation in $\Lambda[\mathbf{x}_\infty]$:
\begin{align*}
    \Bigg\langle
    s_\lambda, ~  \Big(\sum_{\nu \in \Par} s_\nu\Big)^k \Bigg\rangle
    &= \Big\langle
    s_\lambda, ~ \sum_{\bm{\nu} \in (\Par)^k} s_{\nu^{(1)}}\cdots s_{\nu^{(k)}}  \Big\rangle\\
    &=  \Big\langle s_\lambda,~ \sum_\mu \sum_{\bm{\nu} \in (\Par)^k} c_{\bm{\nu}}^\mu s_\mu  \Big\rangle = \sum_{\bm{\nu} \in (\Par)^k} c_{\bm{\nu}}^\lambda.
\end{align*}
Hence, applying the Cauchy identity in $\Lambda[\mathbf{y}_\infty]\otimes\Lambda[\mathbf{z}_\infty]$, we obtain, by Proposition~\ref{p:coefs}:
\begin{align*}
    d_\lambda^\mu(k) = \sum_{\bm{\nu} \in (\Par)^k} c_{\bm{\nu}}^\lambda c_{\bm{\nu}}^\mu &=
    \Bigg\langle
    s_\lambda(\mathbf{y}_\infty)s_\mu(\mathbf{z}_\infty), ~  \Big(\sum_{\nu \in \Par} s_\nu(\mathbf{y}_\infty)s_\nu(\mathbf{z}_\infty)\Big)^k \Bigg\rangle
    \\ &= 
    \Big\langle
    s_\lambda(\mathbf{y}_\infty)s_\mu(\mathbf{z}_\infty), ~  \prod_{i,j} (1-y_jz_i)^{-k} \Big\rangle. \qedhere
\end{align*}
\end{proof}

The third formula for $d_\lambda^\mu(k)$ involves Kronecker products \cite[\S3.2]{Baker}.
\begin{p}\label{p:Kronecker1}
    For all $\lambda,\mu\in\Par$, we have $d_\lambda^\mu(k) = (s_\lambda * s_\mu)(1, \stackrel{k}{\ldots}, 1)$.
\end{p}
\begin{proof}
    Note that $k\mathbf{y} = \mathbf{z}_k\mathbf{y}_{|z_i\mapsto 1}$.
    The generalised Cauchy identity~\eqref{eq:generalised Cauchy} gives
    \[
    s_\lambda[k\mathbf{y}_\infty] = s_\lambda[\mathbf{z}_k\mathbf{y}_\infty]_{|z_i\mapsto 1} =
    \sum_{\mu} (s_\lambda * s_\mu)(\mathbf{z}_k) \cdot s_\mu(\mathbf{y}_\infty)_{|z_i\mapsto 1}.
    \qedhere
    \]
\end{proof}
\begin{cor}\label{cor:d-expresion-ssyt}
    For all $\lambda,\mu\in\Par$, we have $d_\lambda^\mu(k) = \sum_{\nu} g_{\lambda\mu\nu} ~\#\mathrm{SSYT}_k(\nu)$.
\end{cor}
\begin{proof}
    If $|\lambda| \ne |\mu|$ the statement holds trivially. Suppose $|\lambda| = |\mu|$ and expand the Kronecker product to obtain
    \[
    s_\lambda * s_\mu(1, \stackrel{k}{\ldots}, 1) = \sum_\nu g_{\lambda\mu\nu} ~ s_{\nu}(1, \stackrel{k}{\ldots}, 1) 
         = 
        \sum_\nu g_{\lambda\mu\nu} ~\#\mathrm{SSYT}_k(\nu). \qedhere
    \]
\end{proof}

By the hook-content formula, $d_\lambda^\mu(k)$ is therefore a polynomial 
in $k$ \cite[Theorem 5.1]{Brenti}:
    \[
    d_\lambda^\mu(k) 
    = \sum_\nu g_{\lambda\mu\nu} \prod_{u\in Y(\nu)}\frac{k + c(u)}{ h(u)}.
    \]
Crucially, this allows us to define the $k$-fold map for any $k\in\Q_+$; see \S\ref{sec:rational-k-fold}.

We can now deduce necessary conditions for some entries of $D(k)$ to be zero or non-zero.
\begin{cor}
    If $k<\max\left(\frac{\ell(\lambda)}{\ell(\mu)}, \frac{\ell(\mu)}{\ell(\lambda)}\right)$ then $d_\lambda^\mu(k)=0$.
\end{cor}

\begin{proof}
    By the symmetry of $d_\lambda^\mu(k)$ established in Lemma~\ref{lem:symmetries}(1), we may assume $\ell(\lambda)\geq\ell(\mu)$ without loss of generality.
    
    Suppose $k<\ell(\lambda)/\ell(\mu)$ and let $\nu\in\Par$. If $k<\ell(\nu)$, then $\#\text{SSYT}_k(\nu)=0$. If $\ell(\nu)\leq k$, then by the hypothesis, $\ell(\nu)<{\ell(\lambda)}/{\ell(\mu)}$. Recall from \S\ref{sec:preliminaries} that  $g_{\lambda\mu\nu}\neq0$ implies $\ell(\nu)\leq\ell(\lambda)\ell(\mu)$ or, by symmetry of $g_{\lambda\mu\nu}$ in $\lambda,\mu,\nu$, also $\ell(\lambda)\leq\ell(\mu)\ell(\nu)$, which is not satisfied. Hence $g_{\lambda\mu\nu}=0$.

    We have shown that if $k<\ell(\lambda)/\ell(\mu)$ then for every $\nu\in\Par$ the expression $g_{\lambda\mu\nu}\cdot\#\text{SSYT}_k(\nu)$ is zero, and hence so is $d_\lambda^\mu(k)$, by Corollary~\ref{cor:d-expresion-ssyt}.
\end{proof}

\begin{cor}
    If $k\ge|\lambda\cap\mu'|$ then $d_\lambda^\mu(k)\neq0$.
\end{cor}

\begin{proof}
    Recall from \S\ref{sec:preliminaries} that $\max\{\ell(\nu): g_{\lambda\mu\nu}\neq 0\}=|\lambda\cap\mu'|$. Therefore, there exists $\nu_0\in\Par$ such that both $g_{\lambda\mu\nu_0}\neq 0$ and $\ell(\nu_0)=|\lambda\cap\mu'|\le k$, and consequently $\#\text{SSYT}_k(\nu_0)>0$.
    Thus, by Corollary~\ref{cor:d-expresion-ssyt},
    \[
    d_\lambda^\mu(k)=\sum_\nu g_{\lambda\mu\nu}\cdot\#\text{SSYT}_k(\nu)\geq g_{\lambda\mu\nu_0}\cdot \#\text{SSYT}_k(\nu_0)>0. \qedhere
    \]
\end{proof}

\begin{example}\label{eg:D(2)}
    The following is the top-left corner of the infinite matrix $D(2).$
    \[
    \ytableausetup{boxsize=.2em,centertableaux}
    \begin{tikzpicture} \scriptsize
        \node[matrix, row sep={1.2em,between origins}, column sep={1.5em,between origins}] (*) at (0,0)
  {
     & \node{$\varnothing$}; &     \node{\ydiagram{1}}; &\node{\ydiagram{2}}; &\node{\ydiagram{1,1}}; &\node{\ydiagram{3}}; &\node{\ydiagram{2,1}}; &\node{\ydiagram{1,1,1}}; &\node{\ydiagram{4}}; &\node{\ydiagram{3,1}}; &\node{\ydiagram{2,2}}; &\node{\ydiagram{2,1,1}}; &\node{\ydiagram{1,1,1,1}}; \\
\node{$\varnothing$};          & \node{1};&\node{.};&\node{.};&\node{.};&\node{.};&\node{.};&\node{.};&\node{.};&\node{.};&\node{.};&\node{.}; \\
\node{\ydiagram{1}}; & \node{.};&\node{2};&\node{.};&\node{.};&\node{.};&\node{.};&\node{.};&\node{.};&\node{.};&\node{.};&\node{.};&\node{.}; \\
\node{\ydiagram{2}}; & \node{.};&\node{.};&\node{3};&\node{1};&\node{.};&\node{.};&\node{.};&\node{.};&\node{.};&\node{.};&\node{.};&\node{.}; \\
\node{\ydiagram{1,1}}; & \node{.};&\node{.};&\node{1};&\node{3};&\node{.};&\node{.};&\node{.};&\node{.};&\node{.};&\node{.};&\node{.};&\node{.}; \\
\node{\ydiagram{3}}; & \node{.};&\node{.};&\node{.};&\node{.};&\node{4};&\node{2};&\node{.};&\node{.};&\node{.};&\node{.};&\node{.};&\node{.}; \\
\node{\ydiagram{2,1}}; & \node{.};&\node{.};&\node{.};&\node{.};&\node{2};&\node{6};&\node{2};&\node{.};&\node{.};&\node{.};&\node{.};&\node{.}; \\
\node{\ydiagram{1,1,1}}; & \node{.};&\node{.};&\node{.};&\node{.};&\node{.};&\node{2};&\node{4};&\node{.};&\node{.};&\node{.};&\node{.};&\node{.}; \\
\node{\ydiagram{4}}; & \node{.};&\node{.};&\node{.};&\node{.};&\node{.};&\node{.};&\node{.};&\node{5};&\node{3};&\node{1};&\node{.};&\node{.}; \\
\node{\ydiagram{3,1}}; & \node{.};&\node{.};&\node{.};&\node{.};&\node{.};&\node{.};&\node{.};&\node{3};&\node{9};&\node{3};&\node{4};&\node{.}; \\
\node{\ydiagram{2,2}}; & \node{.};&\node{.};&\node{.};&\node{.};&\node{.};&\node{.};&\node{.};&\node{1};&\node{3};&\node{6};&\node{3};&\node{1}; \\
\node{\ydiagram{2,1,1}}; & \node{.};&\node{.};&\node{.};&\node{.};&\node{.};&\node{.};&\node{.};&\node{.};&\node{4};&\node{3};&\node{9};&\node{3}; \\
\node{\ydiagram{1,1,1,1}}; & \node{.};&\node{.};&\node{.};&\node{.};&\node{.};&\node{.};&\node{.};&\node{.};&\node{.};&\node{1};&\node{3};&\node{5}; \\
  };
    \end{tikzpicture}
    \]
\end{example}

\subsection{The monomial basis matrix}
We now turn to the study of the coefficients of the $k$-fold map in the monomial basis. In matrix notation, let $K_{\lambda,\mu} = \langle s_\lambda, m_\mu\rangle=\langle h_\mu, s_\lambda\rangle$ denote the \emph{Kostka coefficient} and let $K = (K_{\lambda,\mu})_{\lambda,\mu\in\Par}$ be the \emph{Kostka matrix}. We study the block-diagonal matrix $B(k) = K^{-1} \cdot D(k) \cdot K$.
\begin{de}
    For $\lambda, \mu \in \Par$, define
    $b_\lambda^\mu(k)$ via $m_\lambda[k \mathbf{y}_\infty] = \sum_\mu b_\lambda^\mu(k) m_\mu(\mathbf{y}_\infty)$.
\end{de}
We give an enumerative formula for $b_\lambda^\mu(k)$, which is dual to~\eqref{eq:alpha in coordinates}. A different enumerative formula in terms of a class of brick tableaux can be found in~\cite[Theorem~4.1]{Brenti}.
\begin{de}
    Given a multiset $A = \{\!\{a_1, a_2,\ldots,a_n\}\!\}$, 
    a \emph{multiset partition} of $A$ is a tuple of multisets (\emph{parts}) $A^{(1)}, A^{(2)}, \ldots, A^{(\ell)}$ such that the union (with multiplicity) of the parts is $A$.

    We define the \emph{weight} of a multiset partition $A^{(1)}, A^{(2)}, \ldots, A^{(\ell)}$ as the tuple $w=(w_1, w_2, \ldots, w_\ell)$, where $w_i$ is the sum of elements of $A^{(i)}$ (with multiplicity). 
    When $k \mid n$, let $\mathcal{P}_{k}(A,w)$ be the set of
    multiset partitions of $A$ of weight $w$ whose parts each have $k$ elements (with multiplicity).
\end{de}
\begin{p}\label{p:m to m}
    Let $|\lambda|=|\mu|$ and $A = \{\!\{\lambda_i  : i = 1, \ldots, k\ell(\mu)\}\!\}$. Then
    \[
    b_{\lambda}^\mu(k) =
    \sum_{\mathcal{P}_{k}(A,\mu)}~
    \prod_{i=1}^{\ell(\mu)}
    \frac{k!}{m_0(A^{(i)})!~m_1(A^{(i)})!\cdots},
    \]
    where $m_j(A^{(i)})$ is the multiplicity of $j$ in $A^{(i)}$.
\end{p}
\begin{proof}
    Let $\ell=\ell(\mu)$. The monomial $y_1^{\mu_1}y_2^{\mu_2}\cdots y_\ell^{\mu_\ell}$ appears in $m_\lambda[k\mathbf{y}_\infty]$ if and only if there is a monomial in $m_\lambda(\mathbf{x}_\infty)$ of the form
    \[
    \Big(x_1^{a_1} x_2^{a_2} \cdots x_{k}^{a_k}\Big)
    \Big(x_{k+1}^{a_{k+1}} \cdots x_{2k}^{a_{2k}}\Big)
    \cdots 
    \Big(x_{(\ell-1)k+1}^{a_{(\ell-1)k+1}} \cdots x_{k\ell}^{a_{k\ell}}\Big)
    \]
    with $a_{(i-1)k+1}+\cdots+a_{ik} = \mu_i$ for all $i$. Hence $A^{(i)} = \{\!\{a_{(i-1)k+1},\ldots,a_{ik}\}\!\}$ defines a multiset partition of $A$ of weight $\mu$ and with $k$ elements in each part.

    For the $i$th parenthesis of the above expression, there are
    \[
    \frac{k!}{m_0(A^{(i)})!~m_1(A^{(i)})!\cdots}
    \]
    ways of reordering the exponents $a_{(i-1)k+1}, \ldots, a_{ik}$.
\end{proof}
In particular, we show shortly that $B(k)$ is lower-triangular. The diagonal is thus $\mathrm{diag}(k^{\ell(\lambda)})_{\lambda\in\Par}$, since $B(k)$ is conjugate to the matrix $X^{-1}P(k)X$ of the $k$-fold map in the power sum basis.
\begin{cor}
    If $b_\lambda^\mu(k) \ne 0$ then $\mu \ge \lambda$ in the lexicographic order.
\end{cor}
\begin{proof}
    If there exists a multiset partition of $\{\!\{\lambda_i  : i = 1, \ldots, k\ell(\mu)\}\!\}$ of weight $\mu$, then clearly $\mu \ge \lambda$.
\end{proof}

Again, we can obtain a formula in terms of Kronecker products~\cite[Proposition 3.2]{Brenti}.
\begin{p}\label{p:b-expresion-Kronecker}
    For all $\lambda,\mu\in\Par$, we have $b_\lambda^\mu(k) = (m_\lambda*h_\mu)(1,\stackrel{k}{\ldots},1)$.
\end{p}
\begin{proof}
It follows from Proposition \ref{p:Kronecker1} and bilinearity that
    \begin{align*}
        b_\lambda^\mu(k) &= \langle m_\lambda[k \mathbf{y}_\infty], h_\mu(\mathbf{y}_\infty) \rangle\\
        &= \sum_{\rho,\nu} K^{-1}_{\lambda,\rho} K_{\nu,\mu} \langle s_\rho[k \mathbf{y}_\infty], s_\nu(\mathbf{y}_\infty) \rangle\\
              &= \sum_{\rho,\nu}K^{-1}_{\lambda,\rho} K_{\nu,\mu} ~(s_\rho*s_\nu)(1, \stackrel{k}{\ldots}, 1)
              = m_\lambda*h_\mu(1, \stackrel{k}{\ldots}, 1). \qedhere
    \end{align*}
\end{proof}

\begin{example}\label{eg:B(2)}
    The following is the top-left corner of the infinite matrix $B(2).$
    \[
    \ytableausetup{boxsize=.2em,centertableaux}
    \begin{tikzpicture} \scriptsize
        \node[matrix, row sep={1.2em,between origins}, column sep={1.5em,between origins}] (*) at (0,0)
  {
     & \node{$\varnothing$}; &     \node{\ydiagram{1}}; &\node{\ydiagram{2}}; &\node{\ydiagram{1,1}}; &\node{\ydiagram{3}}; &\node{\ydiagram{2,1}}; &\node{\ydiagram{1,1,1}}; &\node{\ydiagram{4}}; &\node{\ydiagram{3,1}}; &\node{\ydiagram{2,2}}; &\node{\ydiagram{2,1,1}}; &\node{\ydiagram{1,1,1,1}}; \\
\node{$\varnothing$};          & \node{1};&\node{.};&\node{.};&\node{.};&\node{.};&\node{.};&\node{.};&\node{.};&\node{.};&\node{.};&\node{.};&\node{.};\\
\node{\ydiagram{1}}; & \node{.};&\node{2};&\node{.};&\node{.};&\node{.};&\node{.};&\node{.};&\node{.};&\node{.};&\node{.};&\node{.};&\node{.};\\
\node{\ydiagram{2}}; & \node{.};&\node{.};&\node{2};&\node{.};&\node{.};&\node{.};&\node{.};&\node{.};&\node{.};&\node{.};&\node{.};&\node{.};\\
\node{\ydiagram{1,1}}; & \node{.};&\node{.};&\node{1};&\node{4};&\node{.};&\node{.};&\node{.};&\node{.};&\node{.};&\node{.};&\node{.};&\node{.};\\
\node{\ydiagram{3}}; & \node{.};&\node{.};&\node{.};&\node{.};&\node{2};&\node{.};&\node{.};&\node{.};&\node{.};&\node{.};&\node{.};&\node{.};\\
\node{\ydiagram{2,1}}; & \node{.};&\node{.};&\node{.};&\node{.};&\node{2};&\node{4};&\node{.};&\node{.};&\node{.};&\node{.};&\node{.};&\node{.};\\
\node{\ydiagram{1,1,1}}; & \node{.};&\node{.};&\node{.};&\node{.};&\node{.};&\node{2};&\node{8};&\node{.};&\node{.};&\node{.};&\node{.};&\node{.};\\
\node{\ydiagram{4}}; & \node{.};&\node{.};&\node{.};&\node{.};&\node{.};&\node{.};&\node{.};&\node{2};&\node{.};&\node{.};&\node{.};&\node{.};\\
\node{\ydiagram{3,1}}; & \node{.};&\node{.};&\node{.};&\node{.};&\node{.};&\node{.};&\node{.};&\node{2};&\node{4};&\node{.};&\node{.};&\node{.};\\
\node{\ydiagram{2,2}}; & \node{.};&\node{.};&\node{.};&\node{.};&\node{.};&\node{.};&\node{.};&\node{1};&\node{.};&\node{4};&\node{.};&\node{.};\\
\node{\ydiagram{2,1,1}}; & \node{.};&\node{.};&\node{.};&\node{.};&\node{.};&\node{.};&\node{.};&\node{.};&\node{4};&\node{4};&\node{8};&\node{.};\\
\node{\ydiagram{1,1,1,1}}; & \node{.};&\node{.};&\node{.};&\node{.};&\node{.};&\node{.};&\node{.};&\node{.};&\node{.};&\node{1};&\node{4};&\node{16};\\
  };
    \end{tikzpicture}
    \]
\end{example}

\section{Further directions}\label{further}
There are two natural ways to generalise the results discussed in this article:
\begin{itemize}
    \item Finding an $\SL_2(\C)$ injection that works \emph{without} the divisibility condition.
    \item Extending the construction to an $\SL_N(\C)$ injection.
\end{itemize}
In this section we propose potential candidates for each case, both of which come as generalisations of the $k$-fold map.

\subsection{The \texorpdfstring{$k$}{k}-fold map when \texorpdfstring{$k\not\in\mathbb{Z}$}{k not in Z}}\label{sec:rational-k-fold}
The definition of the $k$-fold map as a plethystic substitution requires $k$ to be a positive integer, while the generalised Foulkes conjecture predicts
\[\Sym_{kd}\Sym^a\C^2\hookrightarrow\Sym_{d}\Sym^{ka}\C^2\]
for any rational $k\ge 1$ such that $ka,kd\in\mathbb{Z}$. Recall from Corollary~\ref{cor:d-expresion-ssyt} that the structure constants of $\kappa_{a,d}^k$ in the Schur basis are polynomials in $k$.

\begin{prob}
    Let $a,d$ be positive integers and $k\in\mathbb{Q}_+$ such that $k\ge 1$ and $ka,kd\in\mathbb{Z}$. Show that the map defined on the Schur basis as
    \begin{align*}
        \kappa_{a,d}^k:\Lambda_{\leq a}[\mathbf{x}_{kd}]&\to\Lambda_{\leq ka}[\mathbf{y}_d]\\
        s_{\lambda}(\mathbf{x}_{kd})&\mapsto\sum_{\mu\in L(d,ka)}\left(\sum_\nu g_{\lambda\mu\nu}\prod_{u\in Y(\nu)}\frac{k+c(u)}{h(u)}\right)s_\mu(\mathbf{y}_d)
    \end{align*}
    for each $\lambda\in L(kd,a)$ is an $\SL_2(\C)$-injection. In particular, that there is an injection 
    \[\Sym_{b}\Sym^a\C^2\hookrightarrow\Sym_{d}\Sym^{c}\C^2\] for arbitrary positive integers $a, b,c,d$ with $ab=cd$ and $a$ being the smallest.
\end{prob}

\subsection{Multisymmetric functions and Foulkes' conjecture}\label{sec:multisym}

The techniques presented in this article may also adapt to the Foulkes' conjecture when $m$ divides $n$ (Conjecture~\ref{conj:0}). We work over an arbitrary field to propose a map, even though Foulkes' conjecture fails over positive characteristic~\cite[Ch.~3]{ODonovan}.\medskip

Foulkes' conjecture is known to hold for $n\le5$ \cite{Thrall,DentSiemons,McKay,CIM}.
Hadamard~\cite{Hadamard} conjectured that for $m\le n$ the following composition is injective, where we have followed our conventions to distinguish dualities. 
\[
\begin{tikzcd}
\Sym_m\Sym_n\CC^N \arrow[d] \arrow["\psi_{m,n}", r, dashed]         & \Sym^n\Sym^m\CC^N                           \\
((\CC^N)^{\otimes n})^{\otimes m} \arrow[r, "\cong"'] & ((\CC^N)^{\otimes m})^{\otimes n} \arrow[u]
\end{tikzcd}
\]
After Howe studied this map in~\cite{Howe}, it became known as the \emph{Howe map}. 
McKay showed that if $\psi_{m,n}$ is injective then $\psi_{m,n+k}$ is injective for all $k\ge0$ \cite{McKay}, which reduced the problem \emph{for each fixed $m$} to a finite computation.
However $\psi_{5,5}$ and $\psi_{6,6}$ are not injective~\cite{MN,CIM}.
\medskip

Consider another possible choice of dualities in order to generalise the $k$-fold map,
\[
\Sym_{kd}\Sym^a\F^N \to \Sym_{d}\Sym^{ka}\F^N.
\]
Following our construction from~\S\ref{sec:preliminaries.poly rings}, we obtain that
\[
\Sym_m\Sym^n\F^N = 
(\F_{\le n}[\mathbf{x}_m^{N-1}])^{S_m}
=: \Lambda^\F_{\le n}[\mathbf{x}_m^{N-1}],
\]
where 
$\mathbf{x}_m^{N-1} = (x_{i,j} \mid 1\le i < N, ~ 1\le j \le m)$ and
$S_m$ acts via $\sigma. x_{i,j} = x_{i,\sigma(j)}$. These polynomials are known in the literature as \emph{multisymmetric polynomials} or MacMahon symmetric polynomials \cite{Rosas,Vaccarino}. 
Treating $\mathbf{x}$ as a matrix of variables $(\mathbf{x}_{i,j})^{1\le i\le N}_{1\le j\le m}$, where $x_{N,j} = 1$ by convention, the action of $A\in\GL_N(\F)$ is given by
\begin{equation}
    \label{eq:actionGLN}
A.f(\mathbf{x}) = 
\prod_{i=1}^{N-1}\big(A'_{*,N}\mathbf{x}_{*,i}\big)^n ~~f\Bigg(\bigg(\frac{A'_{*,j}\mathbf{x}_{*,i}}{A'_{*,N}\mathbf{x}_{*,i}}\bigg)_{i,j}\Bigg),
\end{equation}
where, for a row vector $v$ we let $v'$ denote its transpose, and, for a matrix $M$ we let $M_{*,j}$ denote its $j$th column.

We now break up the alphabet $\mathbf{x}_d^{N-1}$ into a tuple of alphabets
 $(\mathbf{x}_{1,[d]}, \ldots, \mathbf{x}_{N-1,[d]})$, where for each $i=1,\ldots,d$ we have $\mathbf{x}_{i,[d]} = (x_{i,1}, \ldots, x_{i,d})$.
We can define
\[
k\mathbf{x}_d^{N-1} = (k\mathbf{x}_{1,[d]}, \ldots, k\mathbf{x}_{N-1,[d]}).
\]
The $k$-fold map of multisymmetric polynomials is
\begin{align*}
    \kappa^k_{\infty,d} : \Lambda^\F[\mathbf{x}_{kd}^{N-1}] &\to 
    \Lambda^\F[\mathbf{x}_{d}^{N-1}]\\
    f(\mathbf{x}_{kd}^{N-1}) &\mapsto
    f[k\mathbf{x}_{d}^{N-1}]\,.
\end{align*}
Just as before, we denote the restriction to $\Lambda^\F_{\leq a}[\mathbf{x}_{kd}^{N-1}]$ by $\kappa^k_{a,d}$.
\begin{p}
     The image of $\Lambda^\F_{\leq a}[\mathbf{x}_{kd}^{N-1}]$ under $\kappa^k_{a,d}$ lies in $\Lambda^\F_{\leq ka}[\mathbf{x}_{d}^{N-1}]$.	
\end{p}
\begin{proof}
Apply the $k$-fold map to a multisymmetric polynomial $f\in\Lambda^\F_{\leq a}[\mathbf{x}_{kd}^{N-1}]$. 
Fix~$i$.
By Proposition~\ref{p:well-defnd}, the degree of $x_{i,j}$ in each monomial is bounded above by $ka$. This amounts to the claim.
\end{proof}

	\begin{p}
    The $k$-fold map $\kappa_{a,d}^k$ is $\GL_N(\F)$-equivariant.
	\end{p}
	\begin{proof}
        It is a routine check from \eqref{eq:actionGLN} that the $k$-fold map commutes with the action of $\GL_N(\F)$. 
	\end{proof}

\begin{prob}
    Let $\F = \C$ and $a\le d$. Show that the $k$-fold map
    \[
    \Lambda^\C_{\leq a}[\mathbf{x}_{kd}^{N-1}] \to 
    \Lambda^\C_{\leq ka}[\mathbf{x}_{d}^{N-1}]	
    \]
    is injective. In particular, Foulkes' conjecture (Conj.~\ref{conj:0}) holds when $m$ divides~$n$.
\end{prob}

\section*{Acknowledgements}
We thank Abdelmalek Abdesselam for pointing out an error in a previous version of this manuscript, and hinting at a possible relation between $\alpha_d$ and $\kappa_{a,d}^k$. We thank Mark Wildon for his guidance throughout, in particular pointing us to \cite{RSWY} and helping us parse it. We thank Fran\c{c}ois Bergeron for interesting remarks on~\S\ref{sec:matrix}, and Nate Harman for discussions about \S\ref{sec:multisym}.

MG acknowledges financial support from the Austrian Science Fund (FWF) grants 10.55776/P34931 and 10.55776/F1002. ÁG was funded by a University of Bristol Research Training Support Grant. MS gratefully acknowledges financial support from the Heilbronn Institute for Mathematical Research, Bristol.

\bibliographystyle{alpha}
\bibliography{bibliography}

\end{document}